\documentclass[runningheads]{llncs}
\usepackage{amsmath, amssymb}
\usepackage{graphicx}
\usepackage{float}
\usepackage{multirow}
\usepackage{makecell}
\newfloat{lablist}{H}{lol}
\floatname{lablist}{List}

\usepackage{listings}
\usepackage{xcolor}
\usepackage{subcaption}
\usepackage{listings}

\usepackage{algorithm}
\usepackage{algorithmic}
\begin{document}
%
\title{Parallel Cholesky Factorization for Banded Matrices using OpenMP Tasks}
\titlerunning{Parallel Cholesky Factorization for Banded Matrices}
%
\author{Felix Liu \inst{1,2}\orcidID{0000-0001-6865-9379} \and
Albin Fredriksson \inst{2}\and
Stefano Markidis \inst{1}}

\institute{KTH Royal Institute of Technology, Stockholm, Sweden \\
\email{felixliu@kth.se} \and 
RaySearch Laboratories, Stockholm, Sweden}
\maketitle              
\begin{abstract}
  Cholesky factorization is a widely used method for solving linear systems involving symmetric, positive-definite matrices, and can be an attractive choice in applications where a high degree of numerical stability is needed. One such application is numerical optimization, where direct methods for solving linear systems are widely used and often a significant performance bottleneck. An example where this is the case, and the specific type of optimization problem motivating this work, is radiation therapy treatment planning, where numerical optimization is used to create individual treatment plans for patients. To address this bottleneck, we propose a task-based multi-threaded method for Cholesky factorization of banded matrices with medium-sized bands. We implement our algorithm using OpenMP tasks and compare our performance with state-of-the-art libraries such as Intel MKL. Our performance measurements show a performance that is on par or better than Intel MKL (up to \char`\~26\%) for a wide range of matrix bandwidths on two different Intel CPU systems.
\end{abstract}

\keywords{Cholesky factorization, Task-Based Parallelism, Linear Solver}
\section{Introduction}
Cholesky factorization is a well known method for solving linear equations where the matrix is symmetric and positive definite and belongs to a class of algorithms often referred to as \emph{direct} methods for solving linear systems of equations \cite{davis2016survey}. While iterative methods are often considered the state-of-the-art for solving large systems of linear equations, there are still applications where the use of direct methods is the standard, due to issues with ill-conditioning of the linear systems, for instance. Examples of such fields include numerical optimization, where in some algorithms, the systems become increasingly ill-conditioned as the algorithm progresses \cite{gondzio2012interior}. In many applications, matrices involved are not dense but rather have some structure. A common example of structure that arises is banded matrices, where all non-zero elements are located no more than $k$ rows and columns from the main diagonal. In the context of banded matrices, the number $k$ is often referred to as the \emph{bandwidth} of the matrix. With this definition, a diagonal matrix is a banded matrix with bandwidth zero. We emphasize that this use of the term bandwidth is not to be confused with bandwidth when referring to e.g. throughput of memory channels in computer hardware.

The specific problem motivating this work is numerical optimization for radiation therapy treatment planning (see e.g. \cite{unkelbach2015optimization,baumann2016radiation} and references therein for a background on optimization in radiation oncology), where numerical optimization is used to create specialized treatment plans (control parameters for the treatment machine) for individual patients. In such problems, we have observed certain cases which require factorization of banded matrices with bandwidths in the hundreds. Furthermore we have seen this be a significant computational bottleneck, with banded matrix factorization representing more than 50\% of the total time spent in treatment plan optimization.

Of course, optimization problems also arise from a wide range of different application domains, such as operations research, model predictive control among many others. An example from model predictive control is the work of Wang and Boyd \cite{wang2009fast}, where they devise a computational method involving the factorization of matrices with bandwidths up to about 100. We refer readers interested in more details on algorithms for optimization to the review found in \cite{forsgren2002interior}. For a more high-level overview on optimization methods and High-Performance Computing (HPC) we refer the interested reader to the review found in \cite{liu2022survey}.

In this paper, we develop a task-based method for parallelizing Cholesky factorization for banded matrices. We show that our method performs well compared to state of the art libraries on matrices with large bands and further give some discussion and analysis of the performance. We summarize our claimed contributions as follows:
\begin{itemize}
    \item We design and implement a task-based parallel method for Cholesky factorization of banded matrices using OpenMP.
    \item We assess the performance of our method compared to state-of-the-art libraries for matrices with bandwidths between 50-2000, which can be found in optimization problems from radiation therapy.
    \item We demonstrate an up to 26\% performance improvement, on average, compared to state-of-the-art libraries such as Intel MKL.
\end{itemize}
\section{Background} \label{sec:background}
Cholesky factorization is a well-known method for solving systems of linear equations with symmetric, positive-definite matrices. The method works based on the observation that every symmetric positive-definite matrix $A$ admits a factorization of the form $A = LL^T$ where $L$ is a lower triangular matrix. This factorization is called the \emph{Cholesky factorization} and is unique \cite[Ch. 2.7 p. 77]{demmel1997applied}.

One application where banded systems with very large bandwidth can occur is numerical optimization (see e.g. \cite{wright1996applying,wang2009fast}), where many algorithms rely on solving a block-structured linear system of equations in each iteration. The current state-of-the-art in optimization solvers often rely on matrix factorization algorithms to perform this solution step. As an example, in interior point methods \cite{forsgren2002interior} -- a popular choice of algorithm for many constrained optimization problems -- the linear systems to solve will have specific structure depending on structure of the objective and constraints of the optimization problem.. This structure is often exploited in certain cases by, for instance, block-elimination, which may result in the need to solve linear systems with specific structure, where banded structures is one possibility \cite{wang2009fast}.
\subsection{Cholesky Factorization}
\begin{figure}[h!]
\begin{subfigure}[t]{.4\textwidth}
    \centering
    \includegraphics[width=\linewidth]{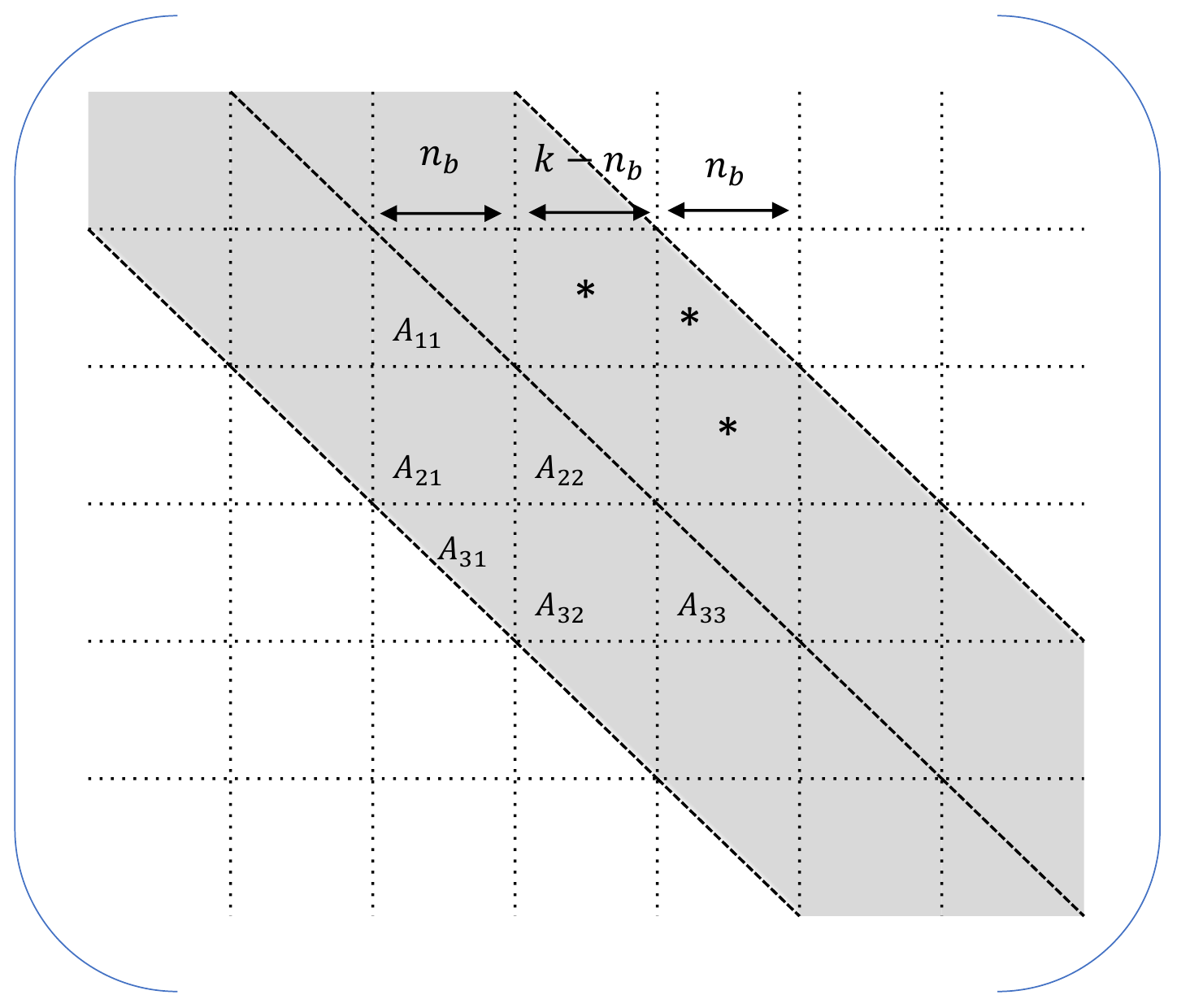}
    \caption{Figure illustrating the block-algorithm for banded matrices from Du Croz et al. \cite{croz1990factorizations}.}
    \label{fig:band_cholesky}
\end{subfigure}
\hfill
\begin{subfigure}[t]{.4\textwidth}
    \centering
    \includegraphics[width=\linewidth]{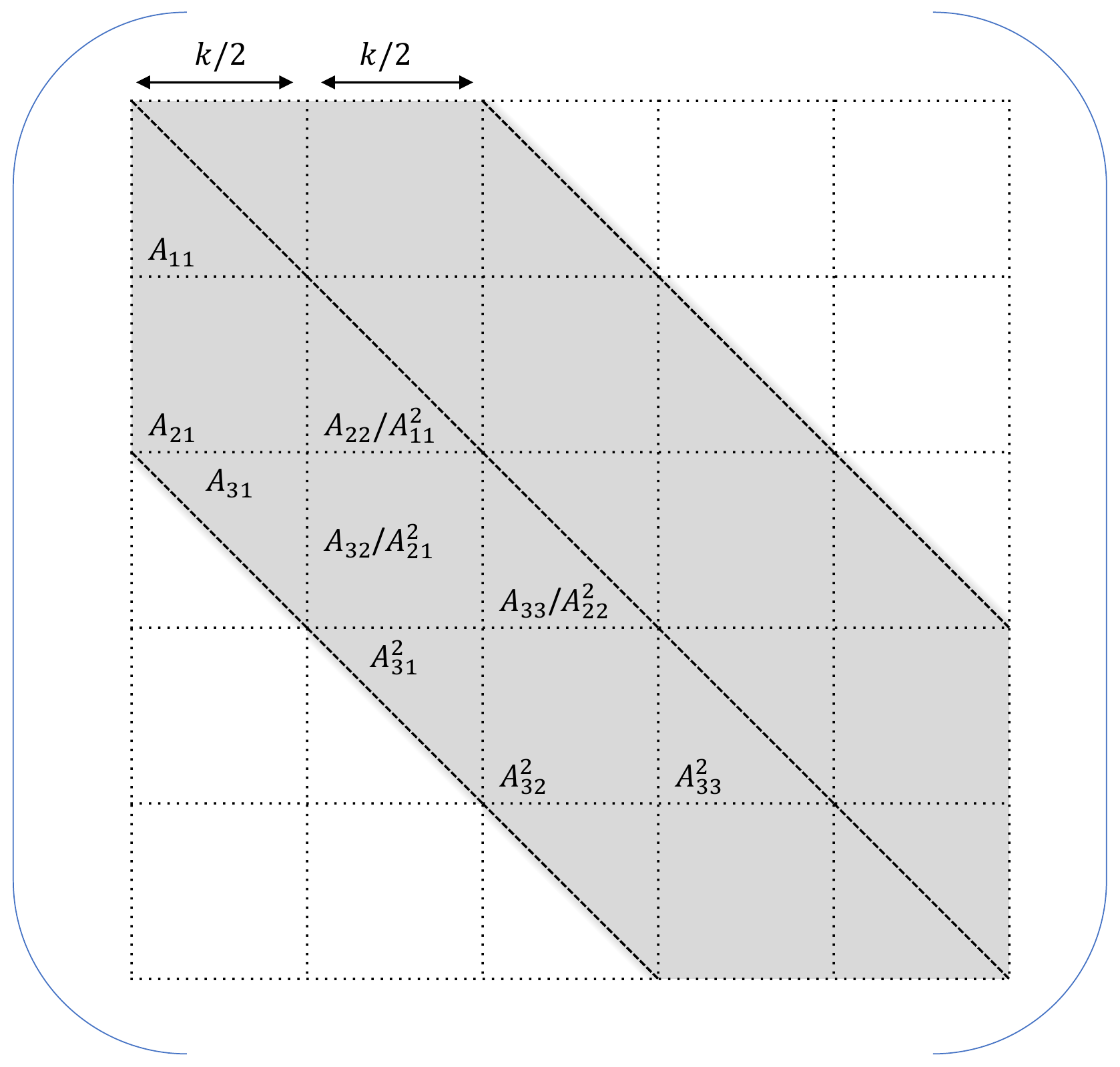}
    \caption{Illustration of two iterations of our parallel algorithm for $3 \times 3$ sub-windows.}
    \label{fig:two_iters}
\end{subfigure}
\caption{Illustration of block algorithms.}
\label{fig:band_cholesky_plots}
\end{figure}

The implementation that provided the initial inspiration for our method is the one proposed by Du Croz et al. in \cite{croz1990factorizations}. The key idea is to organize the computations in the factorization into operations on dense sub-blocks such that the use of level-3 BLAS kernels is enabled. The method of Du Croz et al. divides the current active window in the non-zero band of the matrix into a $3\times3$ \emph{block grid} as illustrated in Fig.~\ref{fig:band_cholesky}. To note is that only the upper triangular part of the $A_{31}$ block lies within the band of the matrix. This is dealt with in the implementation by using an internal square work array with the same dimensions as $A_{31}$, but with the lower triangular part set to zero. A basic outline of the method is:
\begin{enumerate}
    \item \label{en:step_1} Factorize $A_{11}$ into $L_{11}L_{11}^T$ using dense Cholesky
    \item \label{en:step_2} Compute $L_{21} = A_{21}(L_{11}^T)^{-1}$ using DTRSM
    \item Compute $A'_ {22} = A_{22} - L_{21}L_{21}^T$ using DSYRK
    \item Copy the upper triangular part of $A_{31}$ into the square work array
    \item Compute $L_{31} = A_{31}(L_{11}^T)^{-1}$ using DTRSM (store $L_{31}$ in the work array, overwriting $A_{31}$)
    \item Compute $A'_{32} = A_{32} - L_{31}L_{21}^T$ using DGEMM
    \item Compute $A'_{33} = A_{33} - L_{31}L_{31}^T$ using DSYRK
    \item Copy the upper triangular part of $L_{31}$ into the main matrix.
\end{enumerate}
\section{Method} \label{sec:parallel_cholesky}
\noindent
The computational scheme outlined in List~\ref{li:dpbtrf} lends itself naturally to a task-based parallel formulation by inspection of the dependencies between the different steps. For instance, we see that step (2) depends on step (1) for $L_{11}^T$ and step (3) depends on step (2) for $L_{21}$, and so forth. One major drawback of this parallelization scheme is that the $3 \times 3$ active windows may not yield enough parallelism to exploit the hardware to its fullest.

In this work, we extend the $3 \times 3$ block-based algorithm to be able to handle more fine sub-divisions of the active block into $n \times n$ sub-blocks. The benefit of this is a finer task-granularity when parallelizing the algorithm (since each task comprises an operation on a single cell in the grid). One advantage of our method is the use of the standard LAPACK storage format, which lowers the barrier for adoption of our method in existing codes, and avoids possible overhead in converting the matrix to a specialized internal storage format. To note is that we do not use any explicit barriers in our code, but rather all the task dependencies and scheduling is handled by the OpenMP runtime. A pseudo-code implementation of our algorithm is shown in Algorithm~\ref{alg:parallel_cholesky}.
\begin{algorithm}[h!]
\caption{Fine-grained Cholesky factorization}\label{alg:parallel_cholesky}
\begin{algorithmic}[1]
    \REQUIRE $n \geq 0$. \COMMENT{\textit{$n$ is the dimension of the $n \times n$ grid in each iteration}}

\FOR{each active window in the matrix} \label{outer-loop}
    \STATE Factorize $A_{11}$ into $L_{11} L_{11}^T$ using dense Cholesky \label{op1}
    \FOR{$i \gets 2$ to $n - 1$}
        \STATE Compute $L_{i1} = A_{i1} (L_{11}^T)^{-1}$ using DTRSM \label{op2}
            \FOR{$j \gets 2$ to $i-1$}
                \STATE Compute $A_{ij}' = A_{ij} - L_{i1} L_{j1}^T$ using DGEMM \label{op3}
    \ENDFOR
    \STATE Compute $A_{ii}' = A_{ii} - L_{i1} L_{i1}^T$ using DSYRK \label{op4} \\
    \ENDFOR \\
    \COMMENT{\textit{The last row requires special handling, since the bottom-left block is cut off by the banded structure of the matrix.}}
    \STATE{Copy $A_{n1}$ to the work array (a square matrix with bottom left triangle explicitly set to zero).}
    \STATE{Compute $L_{n1} = A_{i1}(L_{11}^T)^{-1}$ using DTRSM (overwriting the value in the work array).}
    \FOR{$k \gets 2$ to $n - 1$}
        \STATE Compute $A_{nk}' = A_{nk} - L_{n1}L_{k1}^T$ using DGEMM (with $L_{n1}$ stored in the work array)
    \ENDFOR
    \STATE Compute $A_{nn}' = A_{nn} - L_{n1}L_{n1}^T$ using DSYRK (with $L_{n1}$ stored in the work array)
\ENDFOR
\end{algorithmic}
\end{algorithm}
\\ \\
The dependency analysis between the different steps shown in Algorithm~\ref{alg:parallel_cholesky} is relatively straightforward when only considering the operations within one \emph{outer iteration} (one iteration of the loop on line~\ref{outer-loop} in Algorithm~\ref{alg:parallel_cholesky}). Each DTRSM operation on the $A_{i1}$ blocks depends on the factorization of the $A_{11}$ block (and the copying to the work array in the case of the bottom left block), the DGEMM operations depend on two of the updates using DTRSM, and the DSYRK operation on the $A_{ii}$ block depends on the DTRSM of the leftmost block on the same block-row. Some special care is required for the bottom block row, due to the use of a temporary array to hold the upper right triangle of the $A_{n1}$ block. To further decrease the amount of synchronization needed, we further extend the analysis of task dependencies to include multiple outer iterations as well. This is made possible by ensuring that the block-grid of the current iteration partially overlaps the block grid from the previous outer iteration, see Fig.~\ref{fig:two_iters} for an illustration for the $n = 3$ case. Thus, the updates on each $A_{ij}$ block that overlaps the previous iteration also depends on the operations on the $A_{i+1,j+1}$ block \emph{from the previous iteration}. This overlapping does however come with the assumption that the dimension $n - 1$ of each $n \times n$ sub-grid evenly divides the bandwidth $k$ in our current implementation.

Finally, the question of how to select the number of sub-blocks to divide each $n \times n$ active window into in the algorithm remains. Recall that we have the requirement that $n - 1$ must divide the bandwidth $k$ of the matrix. As a first prototype, we have implemented a heuristic that tries to select an appropriate value $n$ that fulfills the divisibility requirement while giving the algorithm sufficient parallelism and suitable block sizes to work with. Our heuristic works by selecting a value $n$ that balances the following requirements:
\begin{itemize}
    \item $n - 1$ divides the bandwidth $k$ of the matrix
    \item $n$ is selected such that the block size for the level-3 BLAS operations is approximately 50 by 50.
    \item $n$ is not greater than the number of physical cores of the system.
\end{itemize}
These criteria were selected based on our experimental experience with tuning the performance on our systems. However, this heuristic is still rather crude and may not give optimal performance for all configurations and sizes. As such, users with \emph{a priori} knowledge of the approximate matrix bandwidths and hardware configurations for their use case may tune the number of blocks separately to achieve greater performance. 
\subsection{Implementation Using OpenMP Tasks}
We implement a prototype for our method in C++, where we rely on efficient BLAS implementations to perform the block computations in Algorithm~\ref{alg:parallel_cholesky}. To implement the task-based parallel Cholesky factorization for banded matrices, we use OpenMP tasks with task dependencies, a feature from OpenMP 4.0. To specify data-dependencies between tasks, OpenMP provides \texttt{in}, \texttt{inout} and \texttt{out} clauses to the task pragma, which we use to define the task dependencies described in the previous section. As seen in Listing~\ref{task_snippet}, we use a dummy array to specify the task dependencies, which is not directly accessed by the actual tasks. This works since the OpenMP implementation is agnostic to whether the task dependencies are actually accessed and modified by the tasks, but simply builds the dependency graphs \emph{as if} the data was modified. Thus, the entry at index $(i, j)$ in our dummy array logically represents the block at row index $i$ and column index $j$ in the current active window. For readers interested in other task-based parallel programming frameworks, we refer to the study in \cite{podobas2015comparative}.
\begin{lstlisting}[
    language=C++,
    keywordstyle=\color{blue}\ttfamily,
    stringstyle=\color{red}\ttfamily,
    commentstyle=\color{olive}\ttfamily,
    numbers=left,
    stepnumber=1,
    basicstyle=\small, %or \small or \footnotesize etc.
    caption={Skeletonized C++ code snippet illustrating the implementation of the tasking using OpenMP. Function arguments are omitted for clarity.},
    captionpos=b,
    label=task_snippet,
    columns=fullflexible,
    float=h
]
    char task_dep[BLOCK_DIM][BLOCK_DIM];
    #pragma omp parallel
    #pragma omp single
    {
        for (int i = 0; i < mat_dim; i += nb) {
            #pragma omp task depend(out:task_dep[0][0]) depend(in: task_dep[1][1])
            dpotrf(...);

            for (int blk_i = 1; blk_i < block_dim-1; ++block_i) {
                #pragma omp task depend(in: task_dep[0][0], task_dep[block_i + 1][1]) \
                                 depend(out: task_dep[block_i][0])
                cblas_dtrsm(...);
                for (int block_j = 1; block_j <= block_i; ++block_j) {
                    #pragma omp task depend(in: task_dep[block_i][0], task_dep[block_j][0]) \
                                     depend(out: task_dep[block_i][block_j])
                    cblas_dgemm(...);
                }
                #pragma omp task depend(in: task_dep[block_i][0], \
                                 task_dep[block_i + 1][block_i + 1]) \
                                 depend(out: task_dep[block_i][block_i])
                cblas_dsyrk(...);
            }
        }
    }
\end{lstlisting}
The remainder of our implementation also depends on some BLAS implementation for the level-3 BLAS kernels used in the steps in List~\ref{li:dpbtrf}, as well as an implementation of \texttt{dpotrf} from LAPACK (for step (1)). One important point to note is that our implementation assumes that the bandwidth $k$, defined as the number of super- and sub-diagonals of the matrix, is divisible by $n - 1$, where $n$ is the dimension of the $n \times n$ block grid in each iteration. Our heuristic for selecting $n$ (described in Section~\ref{sec:parallel_cholesky}) selects such a value when possible. Of course, this requirement is impossible to fulfill in certain cases (the bandwidth may for instance be prime). Thus, we have a minimum requirement that the bandwidth of the input matrix is, at least, even such that a division into $3 \times 3$ blocks is valid, which can be ensured by the user by zero-padding their matrix during allocation. Our implementation is available as open source on Github\footnote{https://github.com/felliu/BandCholesky}
\subsection{Performance Model} \label{sec:flop_model}
The algorithm for Cholesky factorization for banded matrices discussed in this paper mainly comprises calls to level-3 BLAS kernels, which are typically highly compute-bound operations. A simple way to judge the performance of a compute-bound algorithm's implementation is to consider the number of floating point operations required to factorize a matrix with given dimensions and (matrix) bandwidth. Since Cholesky factorization is stable without pivoting, the minimum number of floating point operations required to factorize any symmetric positive definite matrix with given dimensions is constant and can be computed relatively exactly. The number of floating point operations required can be computed by considering a left-looking algorithm for Cholesky decomposition, where some loops are truncated from the dense algorithm by the banded structure. Pseudo-code for the algorithm is shown in Algorithm~\ref{alg:cholesky_flops}.
\begin{algorithm}[h!]
\caption{Left-looking Cholesky for banded matrices}\label{alg:cholesky_flops}
 \hspace*{\algorithmicindent} \textbf{Input: }Matrix $A$ with dimension $N$ and bandwidth $k$ \\
 \hspace*{\algorithmicindent} \textbf{Output: }Factorized matrix $L$
\begin{algorithmic}[1]
\FOR{$i \gets 1$ to $N$}
    \FOR{$j \gets \max(1, i - k)$ to $i$}
        \STATE{$t \gets 0$}
        \FOR{$l \gets \max(1, i - k)$ to $j$}
            \STATE{$t \gets t + L(i,l) * L(j, l)$}
        \ENDFOR
        \IF{$i == j$}
            \STATE{$L(i, i) \gets \sqrt{A(i,i) - t}$}
        \ELSE
            \STATE{$L(i, j) \gets (A(i, j) - t) / L(j, j)$}
        \ENDIF
    \ENDFOR
\ENDFOR
\end{algorithmic}
\end{algorithm}
Computing the number of floating point operation required for the factorization can be done in a straightforward way by simply replacing the computational statements in the algorithm with their number of floating point operations and summing for the total value. Let $r = \max(1, i-k)$, then the resulting sum is 
\begin{equation} \label{eq:flop_sum}
\sum_{i=1}^N \sum_{j=r}^i \sum_{l=r}^j 2 + \sum_{i=1}^N \sum_{j=r}^i 2,
\end{equation}
if we consider the square root to be a single floating point operation. While this sum can be evaluated exactly on a computer (which is how we derive the exact FLOP-counts used for the benchmarks), one may also get an approximate value on the order of the number of operations required. We have
\begin{align*}
    &\sum_{i=1}^N \sum_{j=r}^i \sum_{l=r}^j 2 + \sum_{i=1}^N \sum_{j=r}^i 2 \approx \sum_{j=i-k}^i \sum_{l=i-k}^j 2N + 2Nk = \\
    &\sum_{j=i-k}^i 2N(j-(i-k)+1) + 2Nk = 
    2N\sum_{l=1}^k l + 2Nk \\
    &\approx Nk^2 + 2Nk =
    \mathcal{O}(Nk^2),
\end{align*}
where the first approximation is disregarding the truncation in $\max(1, i-k)$, and the second approximation (in the second to last step) uses the observation that the first term is an arithmetic progression.
\section{Experimental Setup}
We evaluate our methods on randomly generated positive-definite (which is ensured by making the matrices diagonally dominant) banded matrices, since the values of the entries do not matter for the number of operations required for Cholesky decomposition (so long as the matrix remains symmetric positive-definite). In all of the experiments below, we fix the dimension of the matrices (the number of rows and columns) to 100,000, and vary the matrix bandwidth.
\subsection{Benchmarking Systems}
In the following we list the benchmarking systems used in this work.
\begin{itemize}
    \item \textbf{Coffee Lake Workstation} is a workstation laptop with a six-core Intel Xeon E-2186M (Coffee Lake) CPU, running Ubuntu 22.04 LTS.
    \item \textbf{Kebnekaise} is an HPC cluster at HPC2N in Umeå, Sweden. Kebnekaise with two Intel Xeon E5-2690v4 (Broadwell) per node. The nodes are running Ubuntu 20.04 LTS.
\end{itemize}
For running benchmarks, we use the Google Benchmark\footnote{\url{https://github.com/google/benchmark}} suite, a C++ library providing different utilities for running (micro)benchmarks. We let the benchmark suite decide the number of iterations to run the benchmark (typically around 10) and then we repeat each run 10 times to gather statistics and estimate noise, all done using the built-in functionality in Google Benchmark. In all the plots below, the median time is reported, to exclude influence from outliers affected by system noise and similar.

For our tasking implementation, we use Intel's OpenMP runtime library, linked with code compiled with GCC (Intel's OpenMP runtime library has a compatibility layer with GNU OpenMP symbols). We compiled our code using GCC 11.2.0 and CMake, with the CMake build set to \texttt{Release} mode (implying optimization level \texttt{-O3} for GCC). The following software versions were used in the experiments:
\begin{itemize}
    \item Intel MKL version 2022.1.0 on Coffee Lake Workstation and 2022.2.0 on Kebnekaise
    \item BLIS version 0.9.0
    \item OpenBLAS version 0.3.20
    \item PLASMA version 21.8.29
\end{itemize}
\section{Results} \label{sec:results}
\begin{table}[h!]
    \centering
    \begin{tabular}{|c|c|c|c|c|c|}
        \hline
        Implementation & \multicolumn{3}{|c|}{Kebnekaise} & \multicolumn{2}{|c|}{Coffee Lake Workstation} \\
        \hline
        Problem setup & \makecell{Low BW\\(50-200)} & \makecell{High BW\\(200-2000)\\Full node} & \makecell{High BW\\(200-2000)\\Single socket} & \makecell{Low BW\\(50-200)} & 
        \makecell{High BW\\(200-2000)} \\
        \hline
        Task Parallel + MKL & \textbf{22.554} & \textbf{218.574} & \textbf{269.133} & \textbf{32.158} & \textbf{168.402} \\
        \hline
        Task Parallel + BLIS & 12.224 & 176.513 & 186.601 & 17.9 & 139.82 \\
        \hline
        MKL Mulithread & 17.256 & 161.078 & 249.022 & 25.089 & 162.2 \\
        \hline
        MKL Sequential & 17.852 & - & - & 25.731 & - \\
        \hline
        PLASMA & - & 127.361 & 128.636 & - & 88.705 \\
        \hline
        \multicolumn{6}{|c|}{} \\
        \hline
        \makecell{\textbf{Improvement over} \\ \textbf{best baseline}} &
        \textbf{26.283\%} & \textbf{35.694\%} & \textbf{8.076\%} & \textbf{24.978}\% & \textbf{3.824\%} \\
        \hline
    \end{tabular}
    \vspace{0.5cm}
    \caption{Summary of \emph{average} performance of the different implementations \emph{in GFLOP/s} over the different matrix bandwidth ranges, test systems and implementations. The specific matrix bandwidths benchmarked are the same as in Fig.~\ref{fig:precdog_perf} and \ref{fig:keb_perf}. The best performing implementation is shown in bold, with the improvements over the best performing baseline (which are sequential and multi-threaded MKL and PLASMA), are shown in the bottom row.}
    \label{tab:perf_table}
\end{table}

\begin{figure}[h!]
\begin{subfigure}[t]{.45\textwidth}
    \centering
    \includegraphics[width=\textwidth]{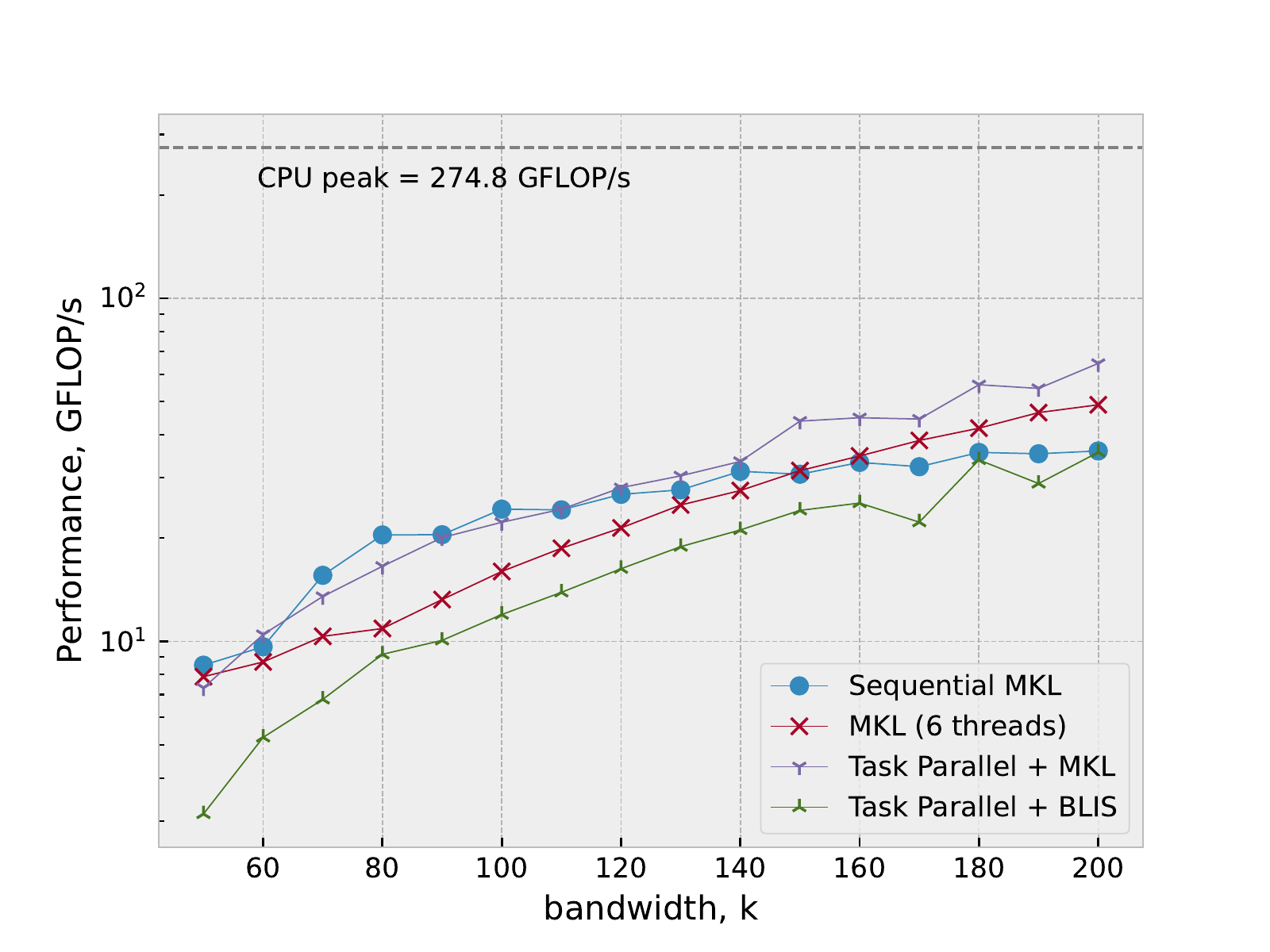}
    \caption{Smaller matrix bandwidths.}
    \label{fig:precdog_lo}
\end{subfigure}
\hfill
    \begin{subfigure}[t]{.45\textwidth}
    \centering
    \includegraphics[width=\textwidth]{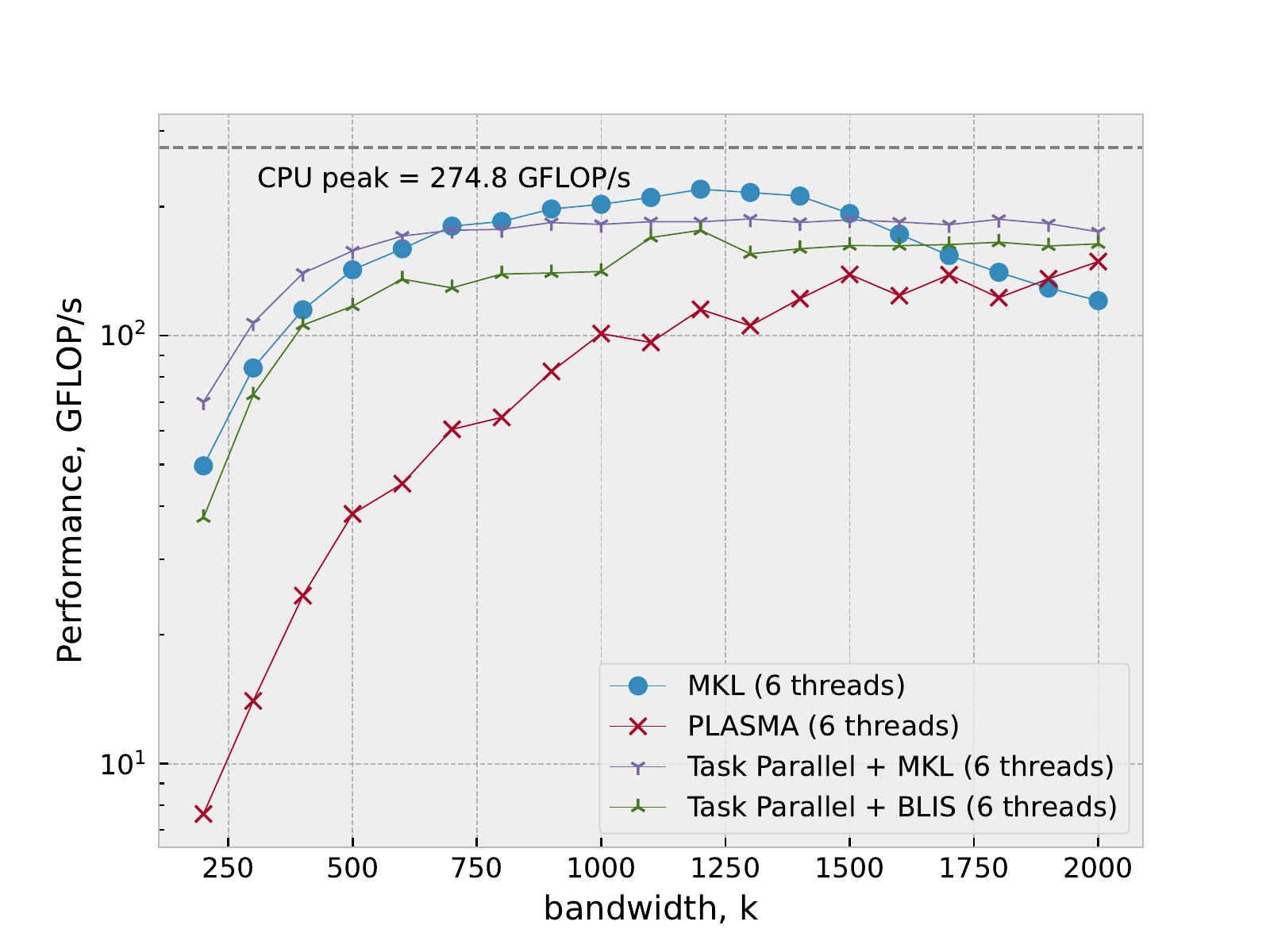}
    \caption{Larger matrix bandwidths.}
    \label{fig:precdog_hi}
\end{subfigure}
\caption{Performance comparison of different Cholesky factorizations on the Coffee Lake Workstation. The performance is shown in GFLOP/s. The left plot shows the performance for smaller matrix bandwidths, and the right for larger matrix bandwidths. Note the log scale on the y-axes.}
\label{fig:precdog_perf}
\end{figure}

\begin{figure}[h!]
\centering
\begin{subfigure}[t]{.45\textwidth}
    \centering
    \includegraphics[width=\textwidth]{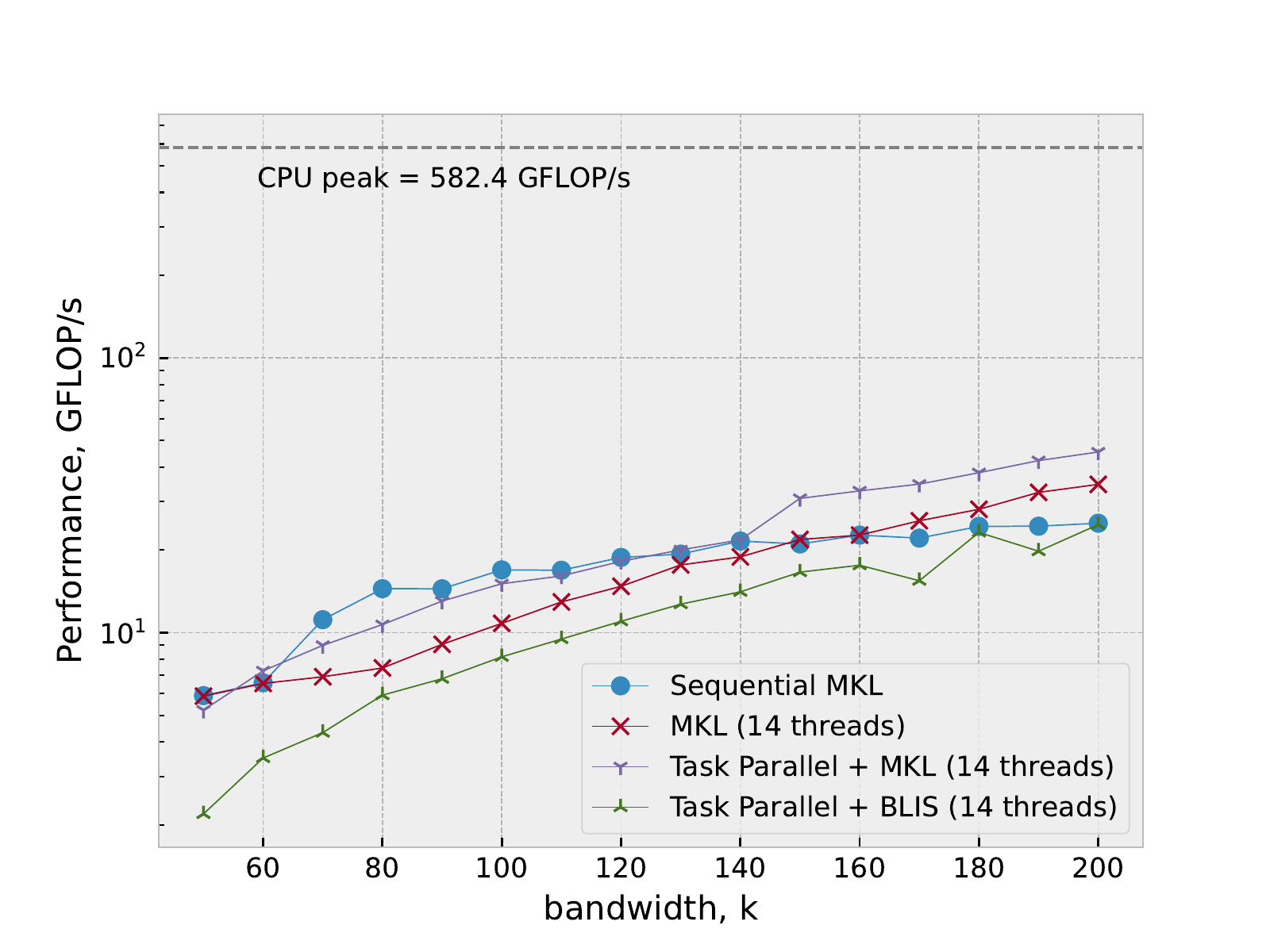}
    \caption{Smaller matrix bandwidths, single CPU socket}
    \label{fig:keb_perf_lo}
\end{subfigure}

\begin{subfigure}[t]{.45\textwidth}
    \centering
    \includegraphics[width=\textwidth]{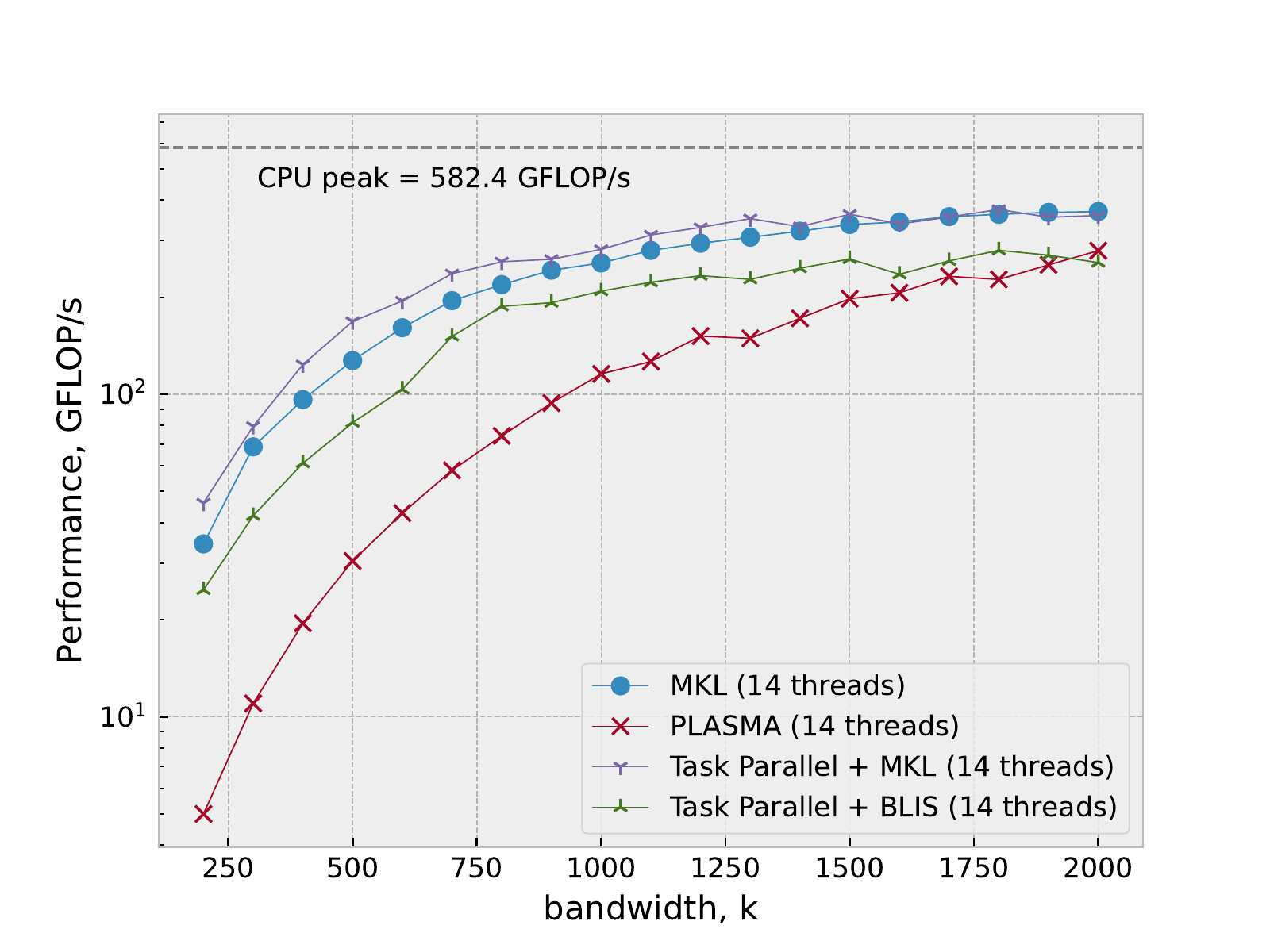}
    \caption{Larger matrix bandwidths, single CPU socket.}
    \label{fig:keb_single_socket}
\end{subfigure}
\hfill
\begin{subfigure}[t]{.45\textwidth}
    \centering
    \includegraphics[width=\textwidth]{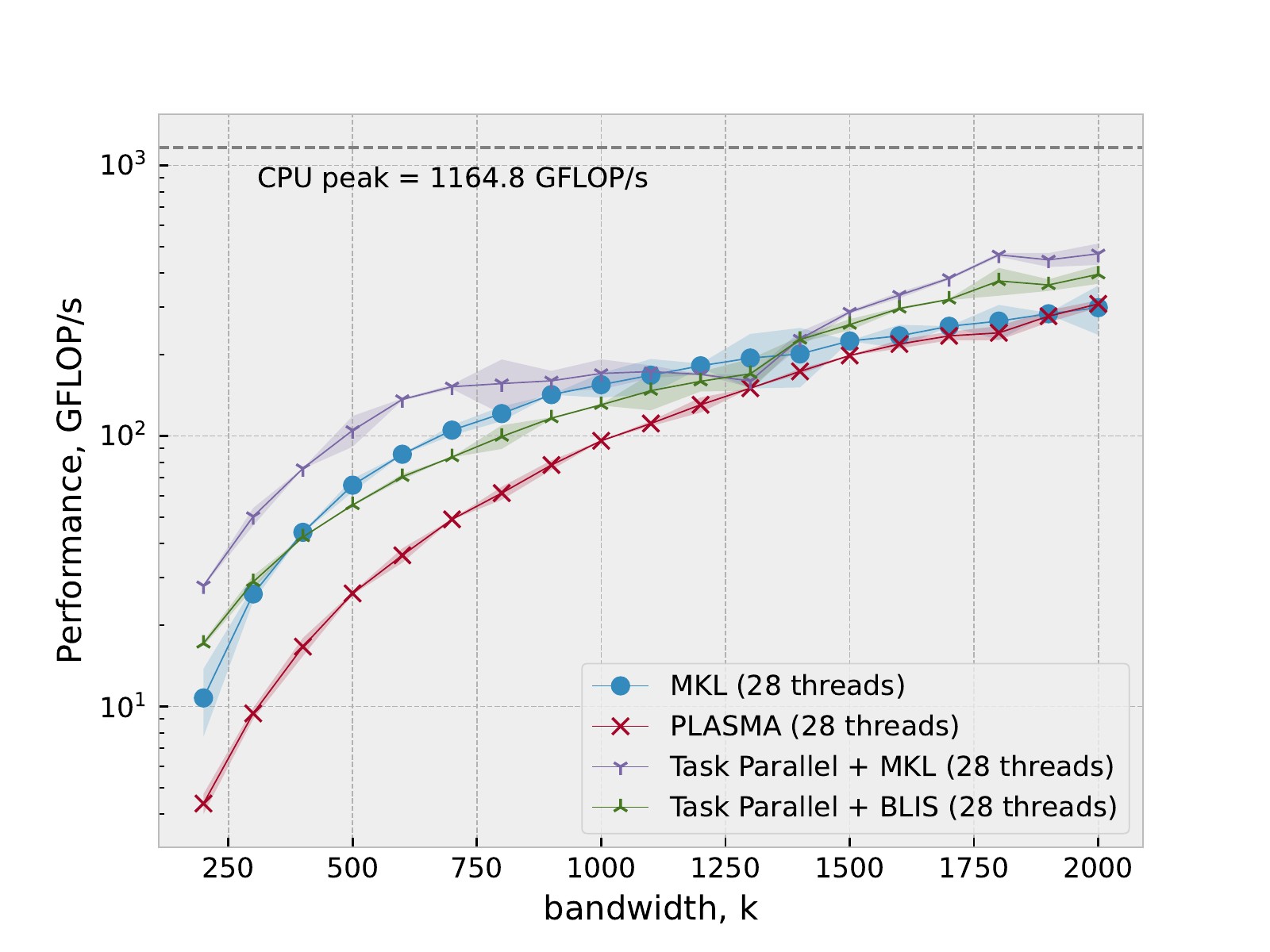}
    \caption{Larger matrix bandwidths, full node.}
    \label{fig:keb_full_node}
\end{subfigure}
\caption{Performance comparison of different Cholesky factorization implementations on Kebnekaise. The performance is shown in GFLOP/s (higher is better). Two types of runs are shown, one using the full node with two Intel Xeon E2690v4 CPUs (28 physical cores) and two corresponding NUMA domains, and two using a single socket (14 physical cores) and a single NUMA domain.}
\label{fig:keb_perf}
\end{figure}
In the following section we present performance results comparing our parallel Cholesky factorization with different state-of-the art libraries (Intel MKL and PLASMA) and settings. Note that our implementation depends on a BLAS implementation as well as an implementation of \texttt{dpotrf} from LAPACK. In some results we will use Intel MKL for these dependencies in our algorithm. These are not to be confused with the stand-alone MKL results, where MKL's \texttt{dpbtrf} (LAPACK kernel for Cholesky factorization of banded matrices) is used. In all of the experiments involving our task-based implementation below, we use the heuristic described in Section~\ref{sec:parallel_cholesky} to decide the dimension $n$ of the block grid in each iteration. In all plots, the performance is shown in GFLOP/s, with the number of floating point operations required computed as described in Section~\ref{sec:flop_model} (in particular, the floating point operation count used to calculate the FLOP/s is the same for all benchmarks). The plots also show the peak performance of the CPU (or node for some Kebnekaise benchmarks) in terms of GFLOP/s in double precision. The values for the peak performance is retrieved from the export compliance metrics provided by Intel for their CPUs (available online\footnote{\url{https://www.intel.com/content/www/us/en/support/articles/000005755/processors.html})}).
\\ \\
The performance plots on the Coffee Lake Workstation are shown in Fig.~\ref{fig:precdog_perf} and the results on Kebnekaise are shown in Fig.~\ref{fig:keb_perf}. On Kebnekaise, the compute nodes have two CPU sockets on different NUMA domains. For the larger matrix bandwidths, we show results using both a single CPU socket (thus avoiding NUMA effects) in Fig.~\ref{fig:keb_single_socket}, and results using the full node in Fig.~\ref{fig:keb_full_node}. In cases where the standard deviation in runtime as reported by Google Benchmark exceeds 5\% (which was only the case for the results in Fig.~\ref{fig:keb_single_socket}), we show a (symmetric) offset of the sample standard deviation as the shaded areas in the plot. The average performance in GFLOP/s over the different matrix bandwidth ranges (50-200 for the low range and 200-2000 for the high range) is summarized in Table~\ref{tab:perf_table}.
\\ \\
We find that our task-based implementation using MKL's BLAS backend is the best performing when considering the average performance across the range of matrix bandwidths, with Intel MKL's \texttt{dpbtrf} being the second best performing in most cases. However, Intel MKL performs better in certain configurations and at certain bandwidths, as we can see in the plots. The difference in performance at different matrix bandwidths for our task-based implementation may be affected by the heuristic used to select sub-block sizes (described in Section~\ref{sec:parallel_cholesky}), which is still rather crude. Furthermore, we see that the performance of our task-based approach using BLIS for the BLAS backend has a rather significant drop in performance compared to using MKL for BLAS. One possible reason for this is that MKL's BLAS level-3 kernels may be better tuned for small matrices (the matrix sizes in each BLAS call will often be approximately $50 \times 50$). This performance difference for smaller matrices has also been observed in previous work \cite{frison2020blas}. PLASMA's performance is lower than MKL in our experiments, which we believe to be caused by overhead in converting the matrix format to PLASMA's internal storage format from the standard LAPACK format used in our benchmarks. On average, we find that our performance improvement relative to Intel MKL is larger for the smaller matrix bandwidth. The average performance across the range of matrix bandwidths is far from the peak performance of the CPUs in all cases, with the size of the gap increasing significantly as the bandwidths of the matrix decreases. For the largest bandwidths, the best performing implementations achieve approximately 70\% of the peak performance of the CPUs.
\section{Related Work} \label{sec:related_work}
Previous research on parallel Cholesky factorization for banded matrices with similar ideas as ours include work by Quintana-Ortí et al. \cite{quintana2008algorithm}, which was implemented in the SuperMatrix \cite{chan2008supermatrix} framework. Our work differs from the work by Quintana-Ortí et al. in that we use the standard packed LAPACK storage format for banded matrices and OpenMP for tasking (OpenMP task implementations were in their infancy at the time the work by Quintana-Ortí et al. was published). We believe this lowers the barrier for adoption in existing codes, and removes the need for potential overhead in converting the matrix into an internal storage format.

The topic of multi-threaded Cholesky factorization is one which has been studied extensively in the literature previously. For example, Remon et al. \cite{remon2006cholesky} performed study on multi-threaded performance for Cholesky factorization of band matrices. In their work, Remon et al. propose some performance optimization by slightly modifying the storage scheme for band matrices used in LAPACK to allow merging some computational steps in List~\ref{li:dpbtrf} into single calls to BLAS-3 kernels, increasing available parallelism for each BLAS kernel invocation. Similar work on small modifications to the LAPACK storage scheme to merge kernels and improve efficiency was also studied by Gustavson et al. \cite{gustavson2008clearer}.

Task parallel banded Cholesky factorization has also been previously studied as part of more extensive efforts to utilize task-based parallelism for different computational kernels in for instance the FLAME \cite{van2009libflame} and PLASMA \cite{dongarra2019plasma} projects. A part of this effort was the previously mentioned work by Quintana-Ortí et al. in \cite{quintana2008algorithm}. On the distributed computing side, parallel algorithms for banded Cholesky has also been studied by e.g. Gupta et al. \cite{gupta1998design}. For general dense matrices the topic has been studied by Dorris et al. \cite{dorris2016task}, where different variants of Cholesky factorization algorithms and their suitability for task-based parallelism were considered. For sparse matrices, task-based parallel Cholesky factorization has been studied in the 1980s by Liu \cite{liu1986computational}, as well as Geist and Ng \cite{geist1989task}, and more recently by Hogg et al. \cite{hogg2010design}.

\section{Conclusions and Future Work} \label{sec:conclusions}
In this paper we have presented our work on evaluating the performance of a task-based parallel algorithm for Cholesky factorization of banded matrices. Our results demonstrated that our method performs, on average, better than state-of-the-art libraries such as Intel MKL for matrices with dimensions and bandwidths similar to those that may arise from our aforementioned optimization problems. However, achieving the optimal performance may depend on a number of factors, including the specific CPU hardware used, the dimensions of the input matrix, among other things.

Furthermore, due to the rapid rise in utilization of GPUs in HPC, investigating the suitability of our algorithm for porting to GPUs is an interesting question. Implementation wise, using task based parallel programming models to target GPUs is possible, for example using the OmpSS programming model \cite{duran2011ompss,bueno2012productive}. Another possible approach is to use the CUDA Graph functionality introduced in CUDA 10, whereby graphs consisting of kernels and their dependencies can be built explicitly and executed on the GPU. Regardless the specific implementation used, one of the main challenges we see is the limited amount of parallelism available in the Cholesky factorization of banded matrices when the size of the bands is modest (our results show that the performance is far from the peak performance of CPUs even at smaller matrix bandwidths). GPUs often require a large amount of available parallelism to run at their peak performance. Thus, the performance benefit of porting to GPUs may be modest, but this is a question we leave for future research.

In conclusion, our implementation performs competitively compared to Intel MKL for our use case, all while keeping a LAPACK-compatible storage scheme for the matrices. Finally, we hope to be able to evaluate our algorithm in a real optimization pipeline for radiation therapy problems in the future to assess the performance improvement in such cases.
\bibliographystyle{splncs04}
\bibliography{References}
\end{document}